\documentclass[10pt, a4paper, twoside]{article}
\usepackage{amssymb}
\usepackage{graphicx}
\usepackage{rotating}

\title{On degenerations between preprojective modules
  over wild quivers}
\author{Roland Olbricht}

\newtheorem{defn}{Definition}[section]
\newtheorem{lem}[defn]{Lemma}
\newtheorem{kor}[defn]{Corollary}
\newtheorem{thm}[defn]{Theorem}

\newcommand{\arsm}[1]{\mathrm{pre}(#1)}

\newcommand{\defect}[2]{d_{#2}(#1)}

\newcommand{\dhom}[2]{\mbox{$\langle #1,#2\rangle$}}
\newcommand{\Ext}{\mathrm{Ext}}
\newcommand{\Hom}{\mathrm{Hom}}

\newcommand{\IN}{\mathbb{N}_0}

\newcommand{\IR}{\mathbb{R}}
\newcommand{\isom}{\cong}
\newcommand{\IZ}{\mathbb{Z}}
\newcommand{\lemref}[3]{\ref{#1}~(#3)}
\newcommand{\mod}{\mbox{mod\ }}
\newcommand{\nullr}[2]{\underline{n}_{#1}(#2)}
\newcommand{\preproj}[1]{\mbox{$\mathcal{P}(\mathcal{#1})$}}
\newcommand{\proj}[1]{P_{#1}}
\newcommand{\pts}[1]{{#1}_0}
\newcommand{\rad}[1]{\mathrm{rad}\ #1}
\newcommand{\sliceQ}[1]{R^q({#1})}
\newcommand{\sliceS}[1]{R^s({#1})}
\newcommand{\spann}[2]{\mbox{$\mathrm{[}#1,#2\mathrm{]}$}}
\newcommand{\reflect}[2]{\sigma_{#1} #2}
\newcommand{\reflectInv}[2]{\sigma^{-1}_{#1} #2}
\newcommand{\supp}{\mbox{supp\ }}
\newcommand{\trde}[1]{\tau^{-1}P_{#1}}
\newcommand{\trdi}[1]{\tau #1}
\newcommand{\trdp}[2]{\tau^{-#1}P_{#2}}
\newcommand{\wild}[2]{\tilde{\tilde{#1}}_{#2}}

\newcounter{vcur}
\newcounter{hcur}
\newcommand{\initcurs}[2]{
  \setcounter{hcur}{#1}
  \setcounter{vcur}{#2}
}

\newcommand{\SBox}{
  \addtocounter{vcur}{-3}
  \put(\value{hcur},\value{vcur}){\line(1,0){6}}
  \put(\value{hcur},\value{vcur}){\line(0,1){6}}
  \addtocounter{vcur}{6}\addtocounter{hcur}{6}
  \put(\value{hcur},\value{vcur}){\line(-1,0){6}}
  \put(\value{hcur},\value{vcur}){\line(0,-1){6}}
  \addtocounter{vcur}{-3}\addtocounter{hcur}{-6}
}
\newcommand{\SDColonSix}{
  \multiput(\value{hcur},\value{vcur})(8,-8){4}{\vector(1,1){8}}
  \multiput(\value{hcur},\value{vcur})(8,-8){3}{\vector(1,-1){8}}
  \put(\value{hcur},\value{vcur}){\vector(1,2){8}}
  \addtocounter{vcur}{8}\addtocounter{hcur}{8}
  \put(\value{hcur},\value{vcur}){\vector(1,-1){8}}
  \addtocounter{vcur}{8}\addtocounter{hcur}{0}
  \put(\value{hcur},\value{vcur}){\vector(1,-2){8}}
  \addtocounter{vcur}{-32}\addtocounter{hcur}{8}
  \put(\value{hcur},\value{vcur}){\vector(1,-2){8}}
  \addtocounter{vcur}{-16}\addtocounter{hcur}{8}
  \put(\value{hcur},\value{vcur}){\vector(1,2){8}}
  \put(\value{hcur},\value{vcur}){\vector(1,-1){8}}
  \addtocounter{vcur}{-8}\addtocounter{hcur}{8}
  \put(\value{hcur},\value{vcur}){\vector(1,1){8}}
  \addtocounter{vcur}{40}\addtocounter{hcur}{-32}
  \addtocounter{vcur}{0}\addtocounter{hcur}{16}
}
\newcommand{\SDColonSixEnde}{
  \multiput(\value{hcur},\value{vcur})(8,-8){4}{\vector(1,1){8}}
  \multiput(\value{hcur},\value{vcur})(8,-8){3}{\vector(1,-1){8}}
  \put(\value{hcur},\value{vcur}){\vector(1,2){8}}
  \addtocounter{vcur}{8}\addtocounter{hcur}{8}
  \put(\value{hcur},\value{vcur}){\vector(1,-1){8}}
  \addtocounter{vcur}{8}\addtocounter{hcur}{0}
  \put(\value{hcur},\value{vcur}){\vector(1,-2){8}}
  \addtocounter{vcur}{-32}\addtocounter{hcur}{8}
  \put(\value{hcur},\value{vcur}){\vector(1,-2){8}}
  \addtocounter{vcur}{-16}\addtocounter{hcur}{8}
  \put(\value{hcur},\value{vcur}){\vector(1,2){8}}
  \addtocounter{vcur}{32}\addtocounter{hcur}{-24}
  \addtocounter{vcur}{0}\addtocounter{hcur}{16}

  \multiput(\value{hcur},\value{vcur})(8,-8){3}{\vector(1,1){8}}
  \multiput(\value{hcur},\value{vcur})(8,-8){2}{\vector(1,-1){8}}
  \put(\value{hcur},\value{vcur}){\vector(1,2){8}}
  \addtocounter{vcur}{8}\addtocounter{hcur}{8}
  \put(\value{hcur},\value{vcur}){\vector(1,-1){8}}
  \addtocounter{vcur}{8}\addtocounter{hcur}{0}
  \put(\value{hcur},\value{vcur}){\vector(1,-2){8}}
  \addtocounter{vcur}{-16}\addtocounter{hcur}{-8}
  \addtocounter{vcur}{0}\addtocounter{hcur}{16}

  \multiput(\value{hcur},\value{vcur})(8,-8){2}{\vector(1,1){8}}
  \put(\value{hcur},\value{vcur}){\vector(1,-1){8}}
  \put(\value{hcur},\value{vcur}){\vector(1,2){8}}
  \addtocounter{vcur}{8}\addtocounter{hcur}{8}
  \put(\value{hcur},\value{vcur}){\vector(1,-1){8}}
  \addtocounter{vcur}{8}\addtocounter{hcur}{0}
  \put(\value{hcur},\value{vcur}){\vector(1,-2){8}}
  \addtocounter{vcur}{-16}\addtocounter{hcur}{-8}
  \addtocounter{vcur}{0}\addtocounter{hcur}{16}
}

\begin{document}

\maketitle

\pagestyle{headings}

\begin{abstract}
  We study minimal degenerations between preprojective modules over wild
  quivers. Asymptotic properties of such degenerations are studied, with
  respect to codimension and numbers of indecomposable direct summands.
  We provide families of minimal disjoint degenerations of arbitrary
  codimension for almost all wild quivers and show that no such examples
  exist in the remaining cases.
\end{abstract}

\section{Introduction}

Given a finite connected quiver $Q$ and a fixed dimension vector, all
representations of $Q$ of this dimension vector are naturally parametrised
by tuples of matrices representing the action of the arrows of $Q$. They
form a vector space, on which a product of general linear groups acts in such
a way that the orbits correspond to the isomorphism classes of
the representations. If the orbit closure of $M$ contains the orbit of $N$
then we say that $M$ degenerates into $N$ or equivalently that $N$ deforms
into $M$. This relation induces an interesting partial order on the set of
isomorphism classes which plays an important role in representation theory.
While there have been substantial advances like
Zwara's module theoretic characterisation of degenerations
\cite{Zwara:1998}, there is still no efficient algorithm available
to compute the closures of the orbit or the degenerations of a module.

In general, a consequence of the degeneration of $M$ into $N$ is that for
any representation $U$ the dimension of the homomorphism space from $U$
to $M$ is at most the dimension of the homomorphism space from $U$ to $N$.
But for certain special classes of representations, called preprojective
representations (to be defined in section 2) we have an equivalence by \cite{Bongartz:1990}. Hence, we can assume without loss of generality
that $M$ and $N$ are disjoint, i.e. have no direct summands in common.

It is shown in \cite{Bongartz:1994} that if $N$ is a minimal
degeneration of $M$, with $M$ and $N$ disjoint, then $N$ is
the direct sum of exactly two indecomposables $U$ and $V$. Here minimal means
that there is no proper degeneration of $M$ of which $N$ is a proper
degeneration. Therefore, to
investigate all minimal disjoint degenerations, it suffices to investigate
all deformations of $U\oplus V$ for any pair of indecomposables $U$ and $V$.

For a quiver of finite representation type this is a finite problem.
It has been solved with the aid of a computer in \cite{Markolf:1990}. In the
case of tame representation type one can still classify the degenerations by
combining a periodicity theorem with computer calculations (see
\cite{Bongartz:2001}). This cannot be expected for wild representation type.
Instead, we turn our attention in this paper to phenomena which occur only for
wild quivers.

Theorem \ref{thm_summands} reveals such a phenomenon: We consider the
preprojective indecomposable modules $U$ and $V$ as vertices in the
preprojective component of the Auslander-Reiten quiver
(see \cite{Auslander_Reiten_Smalo:1995}) as a graph. In
section 2 we introduce a notion of distance for preprojective indecomposable
modules. For any wild quiver and any number $i$ there is a natural number $K(i)$
such that for any pair $U$, $V$ of preprojective indecomposable modules
of a distance greater than $K(i)$ in the Auslander-Reiten quiver
the number of direct summands of each deformation $M$ of $U\oplus V$
is at least $i+1$.

Theorem \ref{thm_codim} is concerned with the codimensions of minimal
disjoint degenerations: while there is a bound of one for quivers of
finite type and a bound of two for quivers of tame type, there exist
minimal disjoint degenerations of arbitrarily high codimension for wild
quivers except for linear quivers with exactly one arrow doubled (denoted
by $K_{m,n}$ in section \ref{Wesentlich}. For those the codimension
is at most two.

The proof of theorem \ref{thm_summands} is brief and
conceptual. It is based on the fact that the kernel of the generalised
Cartan matrix of any wild quiver does not contain any componentwise
nonnegative vector different from zero. In contrast, even the statement
of theorem \ref{thm_codim} already consists of several cases. The proof
requires distinguishing between even more cases which all need different
treatment. But most of the proofs adhere to the same idea so that we only
discuss the case of an extended $\tilde{D}_n$-quiver.

In section 2 we recall some facts about representations and their
deformations. We prove lemma \ref{D.gemischtBeschraenkt} about vector
subspaces in which all vectors have both positive and negative components,
and we formalise a concept already used by Riedtmann in \cite{Riedtmann:1986}
which we denote by deformation shape. Both theorems mentioned above
and the proof of theorem \ref{thm_summands} are the content of section 3,
while the whole of the fourth section is about the proof of
theorem \ref{thm_codim}.

The behaviour of the deformations gives rise to further interesting questions
which are still open: All minimal deformation shapes we find are bounded by
$6$ but there is no obvious reason why they should be bounded at all. Also any
minimal deformation shape seems to be already determined by its support.

\subsection{Acknowledgements}

This paper answers a problem proposed by Klaus~Bongartz and grew out
of a Diploma thesis \cite{Olbricht:2004} written under his supervision.
Markus~Reineke has given me a lot of hints how to improve the readability.

\section{Basic notions and notation}
\label{Begriffe}

Throughout this paper $k$ denotes an algebraically closed field of
arbitrary characteristic. For basic notions and notation we refer to
\cite{Auslander_Reiten_Smalo:1995}. By \emph{quiver} we always mean a
connected directed graph $\mathcal{Q}$ without oriented cycles that
consists of finite sets of vertices $\pts{\mathcal{Q}}$ and arrows
$\mathcal{Q}_1$. This implies that the path algebra $k\mathcal{Q}$
has finite $k$-dimension. We identify representations of the quiver
and $k\mathcal{Q}$-modules and refer to both as \emph{$\mathcal{Q}$-modules}.
The \emph{preprojective indecomposables} are exactly the indecomposable
$\mathcal{Q}$-modules we obtain by applying the transpose of the dual
finitely often to a projective indecomposable $\mathcal{Q}$-module,
or equivalently, all the $\mathcal{Q}$-modules in the
\emph{preprojective component} $\mathcal{Z}$ of the Auslander-Reiten
quiver. On the preprojective indecomposables a partial order is given
by setting $U\preceq V$ if and only if there is a path from $U$ to $V$
inside $\mathcal{Z}$. The length of the shortest path from $U$ to $V$
taken as points in the Auslander-Reiten quiver is called the
\emph{distance} of $U$ and $V$. The set $\spann{U}{V}$ is the set of
all $W$ such that $U\preceq W\preceq V$. In the whole paper we only
consider the full subcategory $\preproj{Q}$ of $\mod k\mathcal{Q}$
which we call \emph{preprojective modules}: These are all finite
direct sums of the preprojective indecomposable $\mathcal{Q}$-modules.
We abbreviate the dual of the transpose of a $\mathcal{Q}$-module $M$
by $\tau M$. Then $\tau^{-1} M$ is the transpose of the dual of the
$\mathcal{Q}$-module $M$. Let $\proj{p}$ be the projective cover of
the simple $\mathcal{Q}$-module associated to a vertex $p$ of $\mathcal{Q}$.

We introduce the abbreviation $\arsm{U}$ for a preprojective
indecomposable $\mathcal{Q}$-module $U$: If $U$ is nonprojective then we set
$\arsm{U}$ equal to the middle term of the Auslander-Reiten sequence
(see \cite{Auslander_Reiten_Smalo:1995}) ending in $U$. If $U$
is projective then we set $\arsm{U}$ equal to $\rad U$.

The matrix
$C\in\IZ^{\pts{\mathcal{Q}}}\times\IZ^{\pts{\mathcal{Q}}}$ defined by
\begin{displaymath}
  C_{pq} := \left\{\begin{array}{cl} 2 & p=q \\
      -n(p,q)-n(q,p) & \textrm{otherwise}
    \end{array}\right.
\end{displaymath}
where we use $n(x,y)$ to refer to the number of arrows from $x$ to $y$ in $\mathcal{Q}$, is called \emph{generalised Cartan matrix of $\mathcal{Q}$}.
For example, if $\mathcal{Q}$ is of type $A$, $D$ or $E$, this is the
Cartan matrix as in \cite{Humphreys:1978}. For any $\mathcal{Q}$
the definition matches the more general one in
\cite{Happel_Preiser_Ringel:1980}.

It is a well-known fact stated e.g. in \cite{Happel_Preiser_Ringel:1980}
that the kernel of $C$ for any wild quiver is a \emph{mixed subspace} of
$\IR^{\pts{\mathcal{Q}}}$, meaning that every non-zero
element of the kernel has positive as well as negative components. We will
need the following remark:

\begin{lem}\label{D.gemischtBeschraenkt}
  Let $U\subset\IR^n$ be a mixed subspace. For every $x\in\IR^n$ the set
  $M=\{x+u|u\in U, x_k+u_k\geq 0\ \textrm{for all}\ k\in\{1,\ldots,n\}\}$
  is bounded.
\end{lem}

\noindent {\em Proof.} \hspace*{6pt}
  Assume that $(v(i))_{i\in\IN}$ is a sequence in $M$ satisfying
  $\lim_{i\to\infty}||v(i)||\to\infty$ where we take $||\cdot||$
  to be the maximum norm of $\IR^n$ for simplicity.
  Set $u(i):=v(i)-x$ and observe that for any component $u(i)_j$ of an
  $u(i)$ the component $-x_j$ is a lower bound. Hence the bounded sequence
  $\frac{1}{||u(i)||}u(i)$ in $U$ has a convergent subsequence which
  converges to a non-zero element of $U$ containing only nonnegative
  components.
\hspace*{12pt}$\Box$

\rule{12pt}{0pt}

For any nonnegative integer $i$ there is a wild quiver whose
Cartan matrix has a kernels of dimension $i$. The situation
is better in the tame case: the kernel of the generalised
Cartan matrix of a tame quiver $\mathcal{Q}$ is always one dimensional
and can be generated by a vector
$\underline{n}_{\mathcal{Q}}\in\IN^{\pts{\mathcal{Q}}}$
which has $1$ as a component. We call $\underline{n}_{\mathcal{Q}}$
the \emph{null-root of $\mathcal{Q}$}. The Coxeter transformation of
$\mathcal{Q}$ does not change $\underline{n}_{\mathcal{Q}}$ and has finite
order, the so-called \emph{Coxeter number of $\mathcal{Q}$}, on the quotient
$\IR^{\pts{\mathcal{Q}}}/\IR\underline{n}_{\mathcal{Q}}$.

A \emph{slice of $\mathcal{Z}$} is a connected full subquiver of
$\mathcal{Z}$ which contains exactly one indecomposable of the $\tau$-orbit of
any projective indecomposable. We introduce some abbreviations for those
types of slices we frequently use: By $\sliceQ{\trdp{i}{j}}$ we denote the
uniquely determined slice with $\trdp{i}{j}$ as its only source. The
slice with the only sink $\trdp{i}{j}$ is called $\sliceS{\trdp{i}{j}}$,
provided it exists. For a function $\delta:\mathcal{Z}\to\IZ$, we abbreviate
by $\delta_R(p)$ the value $\delta(\trdp{i}{p})$ of the uniquely determined
indecomposable of the $\tau$-orbit of $P_p$ that belongs to $R$.

Given a slice $R$ with a source $\trdp{i}{p}$, the slice determined
by the points $(R\setminus\{\trdp{i}{p}\})\cup\{\trdp{(i+1)}{p}\}$ is denoted
by $\reflect{p}{R}$ and is called the \emph{reflection of $R$ at $p$}.

For a quiver $\mathcal{Q}$ and a dimension vector
$d\in\IN^{\pts{\mathcal{Q}}}$ we define the representation variety to be
the affine $k$-variety
\mbox{$\mathcal{D}(\mathcal{Q},d):=\prod_{\alpha\in\mathcal{Q}_1}
  k^{e(\alpha)\times s(\alpha)}$}, where $e(\alpha)$ denotes the ending
vertex of the arrow $\alpha$ and $s(\alpha)$ its starting vertex. The
points of $\mathcal{D}(\mathcal{Q},d)$ correspond to $\mathcal{Q}$-modules
of dimension vector $d$ together with a basis for each of its vector
spaces. The group $G(d):=\prod GL(d_i)$ acts on $\mathcal{D}(\mathcal{Q},d)$
by conjugation such that there is a bijection between the orbits and the
isomorphism classes of $\mathcal{Q}$-modules that have dimension
vector $d$. We observe that the isotropy group of $G(d)$ at each point
$p$ is isomorphic to the group of automorphisms of the corresponding
$\mathcal{Q}$-module $M$. If the orbit of $N$ is contained in the
closure of the orbit of $M$ then we say that $M$ \emph{degenerates}
into $N$ or $N$ \emph{deforms} into $M$ or shortly $M\leq_{deg}N$.
As Riedtmann has observed \cite{Riedtmann:1986} this implies
$\dhom{M}{U}\leq\dhom{N}{U}$ for all $\mathcal{Q}$-modules $U$, where
$\dhom{X}{Y}$ is defined to be the $k$-dimension of $Hom_{\mathcal{Q}}(X,Y)$.
We set the \emph{codimension} of the degeneration of $M$ into $N$
equal to the difference of the dimensions of the orbits of $M$ and
$N$. By a standard dimension calculation we get that
$\dim G(d)p = \dim G(d) - \dhom{M}{M}$, hence $\dhom{N}{N}-\dhom{M}{M}$
is equal to the codimension of the degeneration $M$ in $N$.
We set $M\leq N$ if and
only if $\dhom{M}{U}\leq\dhom{N}{U}$ for all $\mathcal{Q}$-modules $U$.
By a theorem of Auslander \cite{Auslander_Reiten_Smalo:1995},
this is a partial order but in general it differs from the
degeneration order (see \cite{Riedtmann:1986} for an example).

But if $M$ and $N$ are preprojective and have the same dimension vector
it has been proven in \cite{Bongartz:1990} that both orders coincide.
Hence to consider minimal degenerations we can safely assume $M$ and $N$
to have no direct summands in common. Furthermore, it is shown in
\cite{Bongartz:1990} that $N$ can be taken to be $U\oplus V$ for two
indecomposables $U$, $V$ with the property $U\preceq V$. We develop a
notion which is suitable for our calculations:

\begin{defn}\label{B.defshape}
For a quiver $\mathcal{Q}$ and indecomposable preprojective
$\mathcal{Q}$-modules $U$ and $V$
with the property $\tau^{-1} U\preceq V$ we call a function
$\delta:\preproj{Q}\to\IN$ a \emph{deformation shape of $U$ and $V$} if
it satisfies the following conditions:
\begin{itemize}
\item for each $X, Y\in \preproj{Q}$ we have
  $\delta(X\oplus Y)=\delta(X)+\delta(Y)$
\item We have for any indecomposable $\mathcal{Q}$-module $W$: If
  $\delta(W)>0$ then $U\preceq W\preceq \trdi{V}$
\item for any nonprojective indecomposable $\mathcal{Q}$-module $W$ 
  isomorphic neither to $U$ nor to $V$ the \emph{subadditivity inequation}
  $s(\delta,W):=\delta(\arsm{W})-\delta(\trdi{W})-\delta(W)\geq 0$ holds
\item for any projective indecomposable $\mathcal{Q}$-module $W$
  isomorphic neither to $U$ nor to $V$ the \emph{subadditivity inequation}
  $s(\delta,W):=\delta(\arsm{W})-\delta(W)\geq 0$ holds
\item We have $\delta(U)=\delta(\trdi{V})=1$
\end{itemize}
\end{defn}

Let $\delta$ and $\delta'$ be deformation shapes for $U$ and $V$.
We use the term $\delta\leq\delta'$ if and only if $\delta(W)\leq \delta'(W)$
for any indecomposable $W$. This induces a partial order which we refer to by
\emph{deformation order}. A deformation shape $\delta$ is defined to be
\emph{subadditive at $W$} if and only if $s(\delta,W)>0$ and
\emph{strictly additive at $W$} otherwise.

\begin{lem}\label{B.deltaFkt}\ 
\begin{enumerate}
\item\label{B.deltaFkt.add} There is a bijection from the set of
  deformation shapes of $U$ and $V$ to the set of $\mathcal{Q}$-modules $M$ satisfying
  $M< U\oplus V$ which maps a deformation shape $\delta$ onto
  \begin{displaymath}
    M_\delta:=\bigoplus_{W\in\mathcal{Z}, W\not\isom U,V} W^{s(\delta,W)}
  \end{displaymath}
  and whose inverse assigns to a $\mathcal{Q}$-module $M$ the function
  \begin{displaymath}
    \delta_M:=(N\mapsto\dhom{U\oplus V}{N}-\dhom{M}{N})
  \end{displaymath}
\item\label{B.deltaFkt.ord} $M\leq M'$ if and only if $\delta_{M'}\leq \delta_M$
  for $\delta_M$, $\delta_{M'}$ defined as in (a).
\item\label{B.deltaFkt.codim} The codimension of any minimal degeneration of $M$
  into $U\oplus V$ is $\delta_M(M)+1$. It is also called codimension of $\delta_M$.
\end{enumerate}
\end{lem}

\noindent {\em Proof.} \hspace*{6pt}
\begin{enumerate}
\item For any $\delta$ and any indecomposable $W$ isomorphic to
  neither $U$ nor $V$ each $s(\delta, W)$ is nonnegative, hence
  $M_\delta$ is well defined. For each $M$ the deformation shape
  $\delta_M$ is well defined because $M< U\oplus V$ implies
  $\dhom{U\oplus V}{X}\geq\dhom{M}{X}$ for any preprojective $\mathcal{Q}$-module
  $X$ by theorem 3.3 in \cite{Bongartz:1990}.

  By the defining properties of an Auslander-Reiten sequence
  $0\to A\to B\to C\to 0$ we have that for an indecomposable
  $\mathcal{Q}$-module $X$ any non-split morphism $X\to C$ factors
  through $B$. This implies that $\dhom{X}{A}+\dhom{X}{C}-\dhom{X}{B}=0$
  for $X\not\isom C$ and $\dhom{C}{A}+\dhom{C}{C}-\dhom{C}{B}=1$.
  Hence $\dhom{N}{A}+\dhom{N}{C}-\dhom{N}{B}$ is equal to the
  multiplicity of $C$ as a direct summand of $N$. We conclude that
  $s(\delta_M,W)$ is equal to the multiplicity of $W$ as a direct
  summand of $M$, hence $M_{\delta_M}\isom M$.

  To show $\delta_{M_\delta}=\delta$ assume $\delta_{M_\delta}\neq\delta$
  and take a $\preceq$-minimal indecomposable $\mathcal{Q}$-module $X$ with
  $\delta_{M_\delta}(X)\neq\delta(X)$. Then
  \begin{eqnarray*}
    \delta_{M_\delta}(X) & = & \delta_{M_\delta}(\arsm{X})
      - \delta_{M_\delta}(\trdi{X}) - s(\delta_{M_\delta},X) \\
    & = & \delta(\arsm{X}) - \delta(\trdi{X}) - s(\delta,X) = \delta(X)
  \end{eqnarray*}
  yields the required contradiction.

\item $M$ degenerates into $M'$ if and only if
  $\dhom{M}{X}\leq\dhom{M'}{X}$ holds for any preprojective
  $\mathcal{Q}$-module $X$ by theorem 3.3 in \cite{Bongartz:1990}.

\item If $M$ degenerates into $U\oplus V$ minimally then there is an
  exact sequence $0\to U\to M\to V\to 0$ by theorem 4.1 in
  \cite{Bongartz:1990}. Applying $\Hom(U,\_)$ we get
  $\dhom{U}{M}=\dhom{U}{U}+\dhom{U}{V}$ by using that $\Ext^1(U,U)=0$.
  Now (a) enables us to compute
  $\dhom{U\oplus V}{U\oplus V}-\dhom{M}{M}=\dhom{V}{V}+\delta_M(M)$.
\hspace*{12pt}$\Box$
\end{enumerate}

We collect some useful observations about deformation shapes. For this
purpose we need some further notation: Consider a quiver $\mathcal{Q}$
and let $\mathcal{Q}'$ be a full subquiver which is not of finite
representation type. The preprojective component $\mathcal{Z}'$ of the
Auslander-Reiten quiver of $\mathcal{Q}'$ can be embedded naturally
as a full subquiver into the preprojective component $\mathcal{Z}$ of
the Auslander-Reiten quiver of $\mathcal{Q}$.
Given two indecomposables $U'\preceq V'$ in $\mathcal{Z}'$
which are mapped to $U,V\in \mathcal{Z}$ via the embedding and
a deformation shape $\delta'$ of $U'$ and $V'$ we define
$\delta'_!:\mathcal{Z}\to\IN$ to be the extension of $\delta'$ by zero.

\begin{kor}\label{B.zusammenhang}\label{unterkoecher}\label{spiegelungen}
\begin{enumerate}
\item Every deformation shape $\delta$ of $U\oplus V$ has connected support.
\item Assigning to each $\delta'$ the deformation shape $\delta'_!$ yields a
  bijection from the deformation shapes of $U'$ and $V'$ to the deformation
  shapes of $U$ and $V$ which have support only in the image of $\mathcal{Z}'$.
  It preserves the deformation order and the codimension.
\item Let $p\in\pts{\mathcal{Q}}$ be a sink and let $\widetilde{\mathcal{Q}}$ be a
  reflection of $\mathcal{Q}$ at $p$. Then \mbox{$M\leq U\oplus V$} for
  $U,V,M\in \preproj{\widetilde{\mathcal{Q}}}$ is a degeneration if and only if
  the images under the reflection functor in $\preproj{\mathcal{Q}}$ form a
  degeneration. This correspondence preserves the deformation order and codimension.
\end{enumerate}
\end{kor}

\noindent {\em Proof.} \hspace*{6pt}
\begin{enumerate}
\item For any indecomposable $\mathcal{Q}$-module $X$ which is $\preceq$-maximal
  with $\delta(X)>0$, we have $s(\delta,\tau^{-1}X)<0$.

\item Use lemma \ref{B.deltaFkt} for the deformation order and codimension.

\item We assign to each vertex of the preprojective component
  $\widetilde{\mathcal{Z}}$ of $\widetilde{\mathcal{Q}}$ marked by
  $X\in \preproj{\widetilde{\mathcal{Q}}}$ the vertex of $\mathcal{Z}$ marked
  by the image of $X$ under the reflection functor. This yields an embedding of
  $\widetilde{\mathcal{Z}}$ into $\mathcal{Z}$ as a full subquiver. Now proceed
  as above.
\hspace*{12pt}$\Box$
\end{enumerate}

\section{Some special behaviour of wild quivers}
\label{Wesentlich}

We define for any $\mathcal{Q}$-module $M$ its \emph{number of blocks} $\mu(M)$
to be the number of indecomposable nonzero direct summands it can be decomposed
into.

\begin{thm}\label{thm_summands}
  Let $\mathcal{Q}$ be a wild quiver. We define a function $K:\IN\to\IN$
  by setting
  \begin{displaymath}
    K(j) := \min \{\mu(M)|M\leq U\oplus V, U\ \textrm{and}\ V\ 
      \textrm{have distance at least}\ j\}
  \end{displaymath}
  Then $K$ rises monotonically and is unbounded.
\end{thm}

\noindent {\em Proof.} \hspace*{6pt}
  We denote the set of all degenerations between preprojective
  $\mathcal{Q}$-modules $\{(M,U,V)|M\leq U\oplus V\}$ by $D$. Let
  $D_i$ be the set $\{(M,U,V)\in D|\mu(M)\leq i\}$.

  We define a function $t:D\to\IZ^{\pts{\mathcal{Q}}}$ and a function
  $v:D\to\IN^{\pts{\mathcal{Q}}}$ as follows:
  We extend the notation $s(\delta_M,\_)$ used in definition \ref{B.defshape}
  by setting $s(\delta_M, U)=s(\delta_M, V)=-1$. Now for given
  $U$, $V$, $M$ with $M\leq U\oplus V$ we set
  $t((M,U,V))_p := \sum_{l\in\IN} s(\delta_M,\trdp{l}{p})$.
  We define $v((M,U,V))_p := \sum_{l\in\IN} \delta(\trdp{l}{p})$.

  Now $t(D_i)$ is a finite set for any $i$: Any component of any
  $t\in t(D_i)$ is greater or equal than $-2$ and the sum over all
  components $\sum_{p\in\pts{\mathcal{Q}}} t((M,U,V))_p$ is equal to
  $\mu(M)-2$.

  We always get $t((M,U,V))=-Cv((M,U,V))$ with $C$ the Cartan matrix of
  $\mathcal{Q}$: By first applying the definition we obtain for each
  $p\in\pts{\mathcal{Q}}$
  \begin{eqnarray*}
    t((M,U,V))_p\!\! & =\! & \sum_{l\in\IN}\delta_M(\arsm{\trdp{l}{p}})
	- 2\sum_{l\in\IN}\delta_M(\trdp{l}{p}) \\
      & =\! & \delta_M(P_q^{\sum_{q\in\pts{\mathcal{Q}}\setminus\{p\}} n(p,q)}) \\
      & & + \sum_{l\in\IN\setminus\{0\}}\delta_M(
	  \tau^{-l} P_q^{\sum_{q\in\pts{\mathcal{Q}}\setminus\{p\}} n(p,q)}
	  \oplus \tau^{-l+1} P_q^{\sum_{q\in\pts{\mathcal{Q}}\setminus\{p\}} n(q,p)}
	  ) \\
      & & - 2\sum_{l\in\IN}\delta_M(\trdp{l}{p}) \\
      & =\! & \sum_{q\in\pts{\mathcal{Q}}\setminus\{p\}} (n(p,q)+n(q,p))
	\sum_{l\in\IN}\delta_M(\trdp{l}{q})
	- 2\sum_{l\in\IN}\delta_M(\trdp{l}{p})
  \end{eqnarray*}
  We obtain the second equation by using that $\tau$ preserves
  Auslander-Reiten sequences and by using the structure of $\rad P_p$.

  We conclude that $v(D_i)$ is a finite set for any $i$:
  By using that the kernel of the Cartan matrix of a wild quiver
  is always a mixed subspace \cite{Happel_Preiser_Ringel:1980} and lemma
  \ref{D.gemischtBeschraenkt} we obtain that for any given $t$ there
  are only finitely many $v\in\IN^{\pts{\mathcal{Q}}}$ satisfying $t=-Cv$.

  Thus we obtain for each $i$ an upper bound for the possible distance
  of $U$ and $V$ for all $(M,U,V)\in D_i$:
  Given a $w\in\IN^{\pts{\mathcal{Q}}}$, for all $(M,U,V)$ fulfilling
  $v((M,U,V))=w$ the number $\hat{w}:=\sum_{p\in\pts{\mathcal{Q}}} w_p+1$
  is an upper bound for the possible distance of $U$ and $V$.
  Hence $\max\{\hat{w}|w\in v(D_i)\}$ is an upper bound
  for the distance of any $U$ and $V$ such that $(M,U,V)\in D_i$.
\hspace*{12pt}$\Box$

\begin{thm}\label{thm_codim} For any wild quiver $\mathcal{Q}$ we have:
  \begin{enumerate}
  \item\label{codim.kron} There is a global bound on the codimension of any 
    minimal disjoint degeneration of $\mathcal{Q}$ if and only if the
    underlying undirected graph of $\mathcal{Q}$ is one of\\
    \hspace*{20mm}\setlength{\unitlength}{1mm}
    \begin{picture}(72,12)(-24,-4)
      \put(0.5,0){$a_m$}\put(4.5,1){\line(1,0){3.5}}
      \put(8,0){$\ldots$}\put(12,1){\line(1,0){4}}
      \put(16.5,0){$a_1$}\put(20,0.5){\line(1,0){4}}\put(20,1.5){\line(1,0){4}}
      \put(24.5,0){$b_1$}\put(28,1){\line(1,0){4}}
      \put(32,0){$\ldots$}\put(36,1){\line(1,0){4}}
      \put(40.5,0){$b_n$}
      \put(-20,0){$K_{m,n}:=$}
    \end{picture}\\
    for arbitrary $m,n>0$. In this case the bound on the codimension
    is $2$.

  \item\label{codim.standard} Suppose
    $\mathcal{Q}$ contains as a full subquiver 
    a wild quiver in which each pair of vertices is connected by at most
    two arrows and $\mathcal{Q}$ is none of the quivers $K_{m,n}$. Then 
    minimal disjoint degenerations of arbitrary codimension occur in
    $\preproj{\mathcal{Q}}$.
  \end{enumerate}
\end{thm}

\noindent {\em Proof.} \hspace*{6pt}
We first give an overview how the proof proceeds: We have to differentiate
by three different cases. First we treat any quiver $Q$ such that the
underlying undirected graph of $Q$ is a $K_{m,n}$. This proves the ``if'' of (a).

Then we consider quivers whose underlying unoriented graph is of the form\\
\hspace*{40mm}\setlength{\unitlength}{1mm}
\begin{picture}(36,12)(-20,-4)
  \put(0,0){$p$}\put(3,0){\line(1,0){4}}\put(3,1){$\ldots$}
    \put(3,2){\line(1,0){4}}
  \put(8,0){$q$}
  \put(-16,0){$V_m:=$}
\end{picture}\\
consisting of two vertices and $m$ arrows between them. We show that
these admit minimal disjoint degenerations of arbitrarily high comdimension.
For any quiver $Q$ having a $V_m$ as a full subquiver,
corollary~\ref{unterkoecher} then proves that $Q$ admits minimal disjoint
degenerations of arbitrarily high codimension.

Finally we give a list of quivers such that any remaining wild
quiver contains at least one of the quivers from the list as a full
subquiver. Again by corollary~\ref{unterkoecher}, it suffices to give
minimal disjoint degenerations of arbitrary codimension for the
quivers from the list to prove (b) and the ``only if'' of (a).

To show that the codimension of all minimal deformation shapes of $K_{m,n}$
is bounded by $2$ we take a closer look at them: By
corollary~\ref{spiegelungen} without loss of generality, we can fix an
orientation such that $b_1$ is the only sink. First fix
indecomposables $U=P_{a_p}$ and $V=\tau^{-r}P_{a_q}$ such that
$\spann{U}{V}$ has nonempty intersection with the $\tau$-orbits of both
$P_{a_1}$ and $P_{b_1}$. Denote by $W$ the $\preceq$-minimal $\mathcal{Q}$-module
in $\spann{U}{V}$ which is in one of these $\tau$-orbits and by $W'$ the
$\preceq$-maximal $\mathcal{Q}$-module satisfying these conditions. Then the
deformation shape $\delta:\preproj{\mathcal{Q}}\to\IN$ defined by
\begin{displaymath}
  \delta:X\mapsto \left\{
  \begin{array}{cl}
    1 & U\preceq X\preceq W \\
    1 & W'\preceq X\preceq V \\
    1 & X\in\{\trdp{i}{a_1}|i\in\IN\}\cap\spann{U}{V} \\
    1 & X\in\{\trdp{i}{b_1}|i\in\IN\}\cap\spann{U}{V} \\
    0 & \textrm{otherwise}
  \end{array}\right.
\end{displaymath}
is minimal by construction. For any deformation shape $\delta'$ of
$U$ and $V$ we observe that $\delta'(P_{a_{p+j}})\geq \delta'(P_{a_{p+j+1}})$
and by induction using the subadditivity inequation
$\delta'(\tau^k P_{a_j})\geq \delta'(\tau^k P_{a_{j+1}})$ for
$\tau^k P_{a_j}, \tau^k P_{a_{j+1}}\in\spann{U}{\trdi{V}}$. Hence we have
$\delta'(X)\geq 1$ for any $U\preceq X\preceq W$ and we get in a similar
way $\delta'(X)\geq 1$ for any $X$ such that $\delta(X)=1$. It follows
that any deformation shape $\delta'$ of $U$ and $V$ is larger than
$\delta$. The cases of pairs $U=P_{a_p}$, $V=P_{b_q}$ or $U=P_{b_p}$,
$V=P_{a_q}$ or $U=P_{b_p}$, $V=P_{b_q}$ are treated in a similar way.
If $\spann{U}{V}$ has empty intersection with the $\tau$-orbits of
$P_{a_1}$ or $P_{b_1}$ then the subadditivity inequalities admit only
the deformation shape
\begin{displaymath}
  \delta:X\mapsto \left\{
  \begin{array}{cl}
    1 & X\in\spann{U}{\trdi{V}} \\
    0 & \textrm{otherwise}
  \end{array}\right.
\end{displaymath}

On each of the quivers $V_m$ the only minimal deformation shape $\delta$ is
given by $\delta(W)=1$ if and only if $U\preceq W\preceq\trdi{V}$ and
$\delta(W)=0$ everywhere else, due to the connectedness of the support
of $\delta$. When the distance between $U$ and $V$
increases, the codimension of the $\delta$'s increases without bound.

We turn our attention to the list of remaining quivers. Each contains
at least one wild quiver as a full subquiver whose underlying unoriented
graph is one of the following:\\
\setlength{\unitlength}{1mm}
\begin{picture}(37,32)(-21,-16)
  \put(0.5,1){$z$}
  \put(8.5,-7){$a$}\put(4,0){\line(1,-1){4}}
  \put(-7.5,-7){$b$}\put(0,0){\line(-1,-1){4}}
  \put(-7.5,9){$c$}\put(0,4){\line(-1,1){4}}
  \put(8.5,9){$d$}\put(4,4){\line(1,1){4}}
  \put(8.5,1){$e$}\put(4,2){\line(1,0){4}}
  \put(-17.5,1){$S:=$}
\end{picture}
\setlength{\unitlength}{1mm}
\begin{picture}(61,32)(-25,-16)
  \put(0.5,1){$z_1$}\put(4,2){\line(1,0){4}}
  \put(8,1){$\ldots$}\put(12,2){\line(1,0){4}}
  \put(16.5,1){$z_{n-3}$}
  \put(20,0){\line(1,-1){4}}\put(24.5,-7){$d$}
  \put(28,-6){\line(1,0){4}}\put(32.5,-7){$w$}
  \put(20,4){\line(1,1){4}}\put(24.5,9){$v$}
  \put(0,4){\line(-1,1){4}}\put(-7.5,9){$b$}
  \put(0,0){\line(-1,-1){4}}\put(-7.5,-7){$a$}
  \put(-21.5,1){$\wild{D}{n}:=$}
\end{picture}\\
\setlength{\unitlength}{1mm}
\begin{picture}(52,24)(-20,-20)
  \put(1,1){$z$}
  \put(-7.5,1){$b_1$}\put(-4,2){\line(1,0){4}}
  \put(-15,1){$v$}\put(-12,2){\line(1,0){4}}
  \put(0.5,-7){$a_1$}\put(2,-4){\line(0,1){4}}
  \put(0.5,-15){$a_2$}\put(2,-12){\line(0,1){4}}
  \put(8.5,1){$c_1$}\put(4,2){\line(1,0){4}}
  \put(16.5,1){$c_2$}\put(12,2){\line(1,0){4}}
  \put(24.5,1){$w$}\put(20,2){\line(1,0){4}}
  \put(-19.5,-15){$\wild{E}{6}:=$}
\end{picture}
\setlength{\unitlength}{1mm}
\begin{picture}(66,20)(-28,-12)
  \put(1,1){$z$}
  \put(-7.5,1){$b_1$}\put(-4,2){\line(1,0){4}}
  \put(-15.5,1){$b_2$}\put(-12,2){\line(1,0){4}}
  \put(-23,1){$v$}\put(-20,2){\line(1,0){4}}
  \put(0.5,-7){$a_1$}\put(2,-4){\line(0,1){4}}
  \put(8.5,1){$c_1$}\put(4,2){\line(1,0){4}}
  \put(16.5,1){$c_2$}\put(12,2){\line(1,0){4}}
  \put(24.5,1){$c_3$}\put(20,2){\line(1,0){4}}
  \put(32.5,1){$w$}\put(28,2){\line(1,0){4}}
  \put(-27.5,-7){$\wild{E}{7}:=$}
\end{picture}
\setlength{\unitlength}{1mm}
\begin{picture}(84,18)(-28,-12)
  \put(1,1){$z$}
  \put(-7.5,1){$b_1$}\put(-4,2){\line(1,0){4}}
  \put(-15.5,1){$b_2$}\put(-12,2){\line(1,0){4}}
  \put(0.5,-7){$a_1$}\put(2,-4){\line(0,1){4}}
  \put(8.5,1){$c_1$}\put(4,2){\line(1,0){4}}
  \put(16.5,1){$c_2$}\put(12,2){\line(1,0){4}}
  \put(24.5,1){$c_3$}\put(20,2){\line(1,0){4}}
  \put(32.5,1){$c_4$}\put(28,2){\line(1,0){4}}
  \put(41,1){$v$}\put(36,2){\line(1,0){4}}
  \put(48.5,1){$w$}\put(44,2){\line(1,0){4}}
  \put(-27.5,1){$\wild{E}{8}:=$}
\end{picture}\\
\setlength{\unitlength}{1mm}
\begin{picture}(71,16)(-31,-8)
  \put(1,1){$z$}\put(0,2){\line(-1,0){4}}\put(-7,1){$w$}
  \put(4,3){\line(2,1){4}}\put(8.5,5){$b_1$}\put(12,6){\line(1,0){4}}
  \put(16,5){$\ldots$}\put(20,6){\line(1,0){4}}
  \put(24.5,5){$b_{s-1}$}\put(30,5){\line(2,-1){4}}
  \put(4,1){\line(2,-1){4}}\put(8.5,-3){$a_1$}\put(12,-2){\line(1,0){4}}
  \put(16,-3){$\ldots$}\put(20,-2){\line(1,0){4}}
  \put(24.5,-3){$a_{r-1}$}\put(30,-1){\line(2,1){4}}
  \put(35,1){$v$}
  \put(-27,1){$\wild{A}{r,s}:=$}
\end{picture}\\
Except in case $\wild{A}{r,s}$ the reflecting operation on a quiver
is transitive on all possible orientations, so we need to provide
families of minimal disjoint degenerations only for one orientation
for each quiver. In case $\wild{A}{r,s}$ the reflecting operation does
not change the total number of arrows in each direction inside the
circle, so we have to provide families for each number of arrows in
a certain direction.

In all these cases each graph consists of the graph of a tame quiver
and an additional vertex $w$. We will use that observation to give
a family of deformation shapes, prove their minimality and calculate
their codimensions. We will do this explicitly only for the quiver $\wild{D}{n}$
in the following section. For the other quivers in the above list one can
construct such deformation shapes in a similar way.

It remains to give minimal disjoint degenerations for quivers of which
the underlying unoriented graphs are of the form\\
\setlength{\unitlength}{1mm}
\begin{picture}(40,28)(-16,-12)
  \put(1,1){$s$}
  \put(9,9){$a$}\put(4,4){\line(1,1){4}}\put(12,8){\line(1,-1){4}}
  \put(9,1){$b$}\put(4,2){\line(1,0){4}}\put(12,2){\line(1,0){4}}
  \put(9,-7){$c$}\put(4,0){\line(1,-1){4}}\put(12,-4){\line(1,1){4}}
  \put(17,1){$q$}
  \put(-15,1){$AA_5:=$}
\end{picture}\rule{3mm}{0mm}
\setlength{\unitlength}{1mm}
\begin{picture}(40,28)(-16,-12)
  \put(1,1){$s$}\put(4,2){\line(1,0){12}}
  \put(9,9){$a$}\put(4,4){\line(1,1){4}}\put(12,8){\line(1,-1){4}}
  \put(9,-7){$b$}\put(4,0){\line(1,-1){4}}\put(12,-4){\line(1,1){4}}
  \put(17,1){$q$}
  \put(-15,1){$AA_4:=$}
\end{picture}\rule{3mm}{0mm}
\setlength{\unitlength}{1mm}
\begin{picture}(32,28)(-10,-4)
  \put(1,1){$a_1$}\put(4,2){\line(1,0){12}}
  \put(9,9){$a_2$}\put(4,4){\line(1,1){4}}\put(12,8){\line(1,-1){4}}
  \put(17,1){$a_3$}\put(17,4){\line(-1,2){6}}
  \put(9,17){$a_4$}\put(10,16){\line(0,-1){4}}\put(9,16){\line(-1,-2){6}}
  \put(-8,9){$T:=$}
\end{picture}\\
\setlength{\unitlength}{1mm}
\begin{picture}(59,20)(-17,-8)
  \put(0.5,1){$a_1$}
  \put(8.5,1){$a_2$}\put(4,1.5){\line(1,0){4}}\put(4,2.5){\line(1,0){4}}
  \put(16,1){$\ldots$}\put(12,2){\line(1,0){4}}
  \put(24.5,1){$a_{n-1}$}\put(20,2){\line(1,0){4}}
  \put(34.5,1){$a_n$}\put(31,1.5){\line(1,0){3}}\put(31,2.5){\line(1,0){3}}
  \put(-17,1){$KK_n:=$}
\end{picture}\rule{5mm}{0mm}
\setlength{\unitlength}{1mm}
\begin{picture}(55,20)(-17,-8)
  \put(0.5,1){$z_1$}
  \put(8.5,1){$z_2$}\put(4,1.5){\line(1,0){4}}\put(4,2.5){\line(1,0){4}}
  \put(16,1){$\ldots$}\put(12,2){\line(1,0){4}}
  \put(24,1){$z_{n-2}$}\put(20,2){\line(1,0){4}}
  \put(35,5){$a$}\put(31,3.5){\line(2,1){4}}
  \put(35,-3){$b$}\put(31,0,5){\line(2,-1){4}}
  \put(-17,1){$KD_n:=$}
\end{picture}\\
\begin{picture}(59,20)(-17,-8)
  \put(0.5,1){$z_1$}
  \put(8.5,1){$z_2$}\put(4,1.5){\line(1,0){4}}\put(4,2.5){\line(1,0){4}}
  \put(16,1){$\ldots$}\put(12,2){\line(1,0){4}}
  \put(24,1){$z_{n-2}$}\put(20,2){\line(1,0){4}}
  \put(35,5){$a$}\put(31,3.5){\line(2,1){4}}\put(36,0){\line(0,1){4}}
  \put(35,-3){$b$}\put(31,0.5){\line(2,-1){4}}
  \put(-17,1){$KA_n:=$}
\end{picture}\rule{5mm}{0mm}
\setlength{\unitlength}{1mm}
\begin{picture}(55,20)(-17,-8)
  \put(1,1){$z$}
  \put(4,3){\line(2,1){4}}\put(8.5,5){$b_1$}\put(12,6){\line(1,0){4}}
  \put(16,5){$\ldots$}\put(20,6){\line(1,0){4}}
  \put(24.5,5){$b_{q-1}$}\put(30,5){\line(2,-1){4}}
  \put(4,0.5){\line(2,-1){4}}\put(4,1.5){\line(2,-1){4}}
  \put(8.5,-3){$a_1$}\put(12,-2){\line(1,0){4}}
  \put(16,-3){$\ldots$}\put(20,-2){\line(1,0){4}}
  \put(24.5,-3){$a_{p-1}$}\put(30,-1){\line(2,1){4}}
  \put(35,1){$c$}
  \put(-17,1){$AK_{p,q}:=$}
\end{picture}\\
\setlength{\unitlength}{1mm}
\begin{picture}(40,20)(-16,-12)
  \put(1,1){$s$}\put(4,2){\line(1,0){12}}
  \put(9,-7){$b$}\put(4,-0.5){\line(1,-1){4}}\put(12,-4.5){\line(1,1){4}}
  \put(4,0.5){\line(1,-1){4}}\put(12,-3.5){\line(1,1){4}}
  \put(17,1){$q$}
  \put(-15,1){$\tilde{AK}_3:=$}
\end{picture}
\setlength{\unitlength}{1mm}
\begin{picture}(47,20)(-23,-8)
  \put(1,1){$z$}
  \put(4,3){\line(2,1){4}}\put(9,5){$b$}
  \put(12,4.5){\line(2,-1){4}}\put(12,5.5){\line(2,-1){4}}
  \put(4,0.5){\line(2,-1){4}}\put(4,1.5){\line(2,-1){4}}
  \put(9,-3){$a$}\put(12,-1){\line(2,1){4}}
  \put(17,1){$c$}
  \put(-19,1){$\tilde{AK}_4:=$}
\end{picture}\\
These are easier to find. Details can be found in \cite{Olbricht:2004}.
We give only one example here: Let $\mathcal{Q}$ be the quiver with the
unoriented graph $KK_3$ and the only sink $a_2$, let $m$ be a positive
integer, $U:=\proj{a_2}$ and $V:=\trdp{m}{a_2}$. Then
\begin{displaymath}
  \delta:\trdp{i}{j}\mapsto\left\{\begin{array}{cl}
      1 & j=a_1, i\in\{2k|k\in\IN, 2k < m-1)\}\\
      1 & j=a_2, i\in\{0, \ldots, m-1\}\\
      1 & j=a_3, i\in\{2k+1|k\in\IN, 2k+1 < m-1)\}\\
      0 & \textrm{otherwise}
    \end{array}\right.
\end{displaymath}
is a deformation shape of codimension $m$. Its minimality can be proven
immediately by the subadditivity equations.
\hspace*{12pt}$\Box$

\section{Families of minimal degenerations}
\label{Rechnungen}
\label{R.fastZahm}

We provide in this section a family of minimal deformation shapes
of any codimension greater or equal $2$ for quivers whose
underlying unoriented graph is of the form $\wild{D}{n}$.
Due to corollary~\ref{spiegelungen}, such a family for some orientation
can easily be transformed into a family with the same property for
another orientation. Hence we fix an orientation specified by making
$z_1$ the only sink. Let $c$ be the Coxeter number of $\wild{D}{n}$.
Now for given $m$ we set $U:=P_{z_1}$ and $V:=\trdp{((m-1)(c+1)+n-2)}{z_1}$
and define a deformation shape $\delta_m$ of $U$ and $V$.

For the purpose of defining $\delta_m$ and proving its minimality we cover
$\spann{U}{\tau V}$ by segments:
\begin{itemize}
\item the initial segment consisting of $\preceq$-predecessors of any
  indecomposable in $\sliceQ{B^q_0}$
\item the final segment consisting of all $\preceq$-successors of
  $\sliceS{B^s_{m-1}}$
\item the tame segments between $\sliceQ{B^q_i}$ and $\sliceS{B^s_{i+1}}$
  for $0\leq i\leq m-2$ and
\item the wild segments between $\sliceS{B^s_i}$ and $\sliceQ{B^q_i}$
  for $1\leq i\leq m-1$
\end{itemize}
where $B^q_i := \trdp{(i(c+1)+1)}{z_{n-3}}$
and $B^s_i := \trdp{i(c+1)}{z_{n-3}}$. Inside the tame segments the
deformation shapes $\delta_m$ have support only inside the $\tau$-orbits
of $\pts{\mathcal{Q}}\setminus\{w\}$.

Unfortunately, the proof is intricate. So we start with an overview over
the strategy of the proof: The proof will start with the definition of
the defect of a deformation shape at a slice. This definition extends
the usual definition of the defect for tame quivers onto wild quivers
that consist of a tame quiver as a full subquiver and an additional vertex
$w$.

Now the difference between the defects of a deformation shape
on the leftmost slice and the rightmost slice of a segment reveals
essential properties of the deformation shape: We conclude from
lemma~\ref{R.defekt} that on any tame segment
the defect of a deformation shape that has support only inside the
$\tau$-orbits of $\pts{\mathcal{Q}}\setminus\{w\}$ is on
the leftmost slice always greater or equal than on the rightmost slice of
this segment. Furthermore the deformation shape is strictly additive on
the $\tau$-orbits of $\pts{\mathcal{Q}}\setminus\{w\}$ in the
tame segment if and only if these defects are equal.

An essential upper bound for the defect of the leftmost slice of the
final segment can now be computed by induction: Lemma~\ref{dttn}
starts the induction on the initial segment in part (a) and makes
it work on the wild segments in part (b). To complete the induction,
observe that in any tame segment between $\sliceQ{B^q_i}$ and
$\sliceS{B^s_{i+1}}$ the deformation shape $\delta$ is properly
subadditive or it is subadditive at $\trdp{((i+1)(c+1)+1)}{v}$
because $\delta\leq\delta_m$ and $\delta_m(\trdp{((i+1)(c+1)+1)}{v})=0$.
This yields that the defect of $\delta$ on each leftmost slice
of a tame segment is at most $0$.

To proceed, we will use the symmetry inherent to the conditions for
deformation shapes to define for any deformation shape $\delta$ of
$U$ and $V$ the reflection $\overline{\delta}$. This is again a
deformation shape of $U$ and $V$ such that
$\defect{\overline{\delta}}{\overline{R}^{\delta}}=-\defect{\delta}{R}$.

The proof of the minimality is completed by lemma~\ref{dttn}~(c).

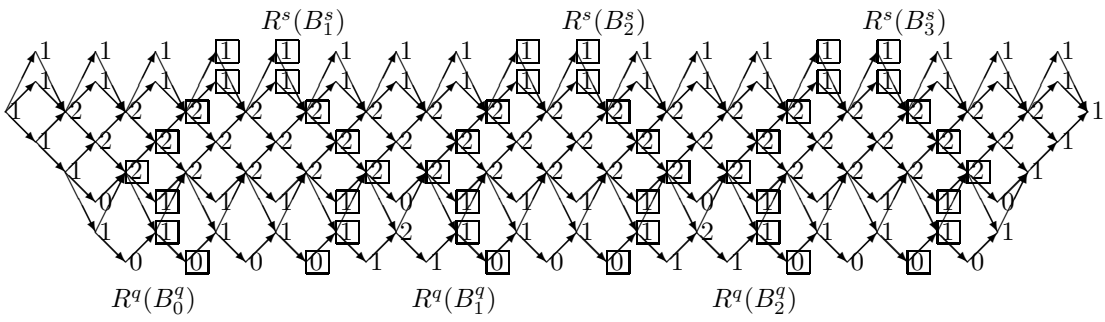
\begin{figure}
\begin{sideways}
  \setlength{\unitlength}{0.5mm}
\begin{picture}(276,8)(0,34)
  \initcurs{0}{0}
  \SDColonSix\SDColonSix\SDColonSix\SDColonSix\SDColonSix
  \SDColonSix\SDColonSix\SDColonSix\SDColonSix\SDColonSix
  \SDColonSix\SDColonSix\SDColonSix\SDColonSix\SDColonSix
  \SDColonSixEnde

  \multiput(9,14)(16,0){18}{1}
  \multiput(9,6)(16,0){18}{1}
  \put(1,-2){1}\multiput(17,-2)(16,0){17}{2}\put(289,-2){1}
  \put(9,-10){1}\multiput(25,-10)(16,0){16}{2}\put(281,-10){1}
  \put(17,-18){1}\multiput(33,-18)(16,0){15}{2}\put(273,-18){1}

  \put(25,-26){0}\multiput(41,-26)(16,0){4}{1}
  \put(105,-26){0}\multiput(121,-26)(16,0){4}{1}
  \put(185,-26){0}\multiput(201,-26)(16,0){4}{1}\put(265,-26){0}

  \multiput(25,-34)(16,0){5}{1}\put(105,-34){2}
  \multiput(121,-34)(16,0){4}{1}\put(185,-34){2}
  \multiput(201,-34)(16,0){5}{1}

  \multiput(33,-42)(16,0){4}{0}
  \multiput(97,-42)(16,0){2}{1}\multiput(129,-42)(16,0){3}{0}
  \multiput(177,-42)(16,0){2}{1}\multiput(209,-42)(16,0){4}{0}

  \initcurs{48}{0}\SBox\initcurs{56}{8}\SBox\initcurs{56}{16}\SBox
  \initcurs{40}{-8}\SBox\initcurs{32}{-16}\SBox\initcurs{40}{-24}\SBox
  \initcurs{40}{-32}\SBox\initcurs{48}{-40}\SBox
  \put(28,-52){$\sliceQ{B^q_0}$}

  \initcurs{80}{0}\SBox\initcurs{72}{8}\SBox\initcurs{72}{16}\SBox
  \initcurs{88}{-8}\SBox\initcurs{96}{-16}\SBox\initcurs{88}{-24}\SBox
  \initcurs{88}{-32}\SBox\initcurs{80}{-40}\SBox
  \put(68,22){$\sliceS{B^s_1}$}

  \initcurs{128}{0}\SBox\initcurs{136}{8}\SBox\initcurs{136}{16}\SBox
  \initcurs{120}{-8}\SBox\initcurs{112}{-16}\SBox\initcurs{120}{-24}\SBox
  \initcurs{120}{-32}\SBox\initcurs{128}{-40}\SBox
  \put(108,-52){$\sliceQ{B^q_1}$}

  \initcurs{160}{0}\SBox\initcurs{152}{8}\SBox\initcurs{152}{16}\SBox
  \initcurs{168}{-8}\SBox\initcurs{176}{-16}\SBox\initcurs{168}{-24}\SBox
  \initcurs{168}{-32}\SBox\initcurs{160}{-40}\SBox
  \put(148,22){$\sliceS{B^s_2}$}

  \initcurs{208}{0}\SBox\initcurs{216}{8}\SBox\initcurs{216}{16}\SBox
  \initcurs{200}{-8}\SBox\initcurs{192}{-16}\SBox\initcurs{200}{-24}\SBox
  \initcurs{200}{-32}\SBox\initcurs{208}{-40}\SBox
  \put(188,-52){$\sliceQ{B^q_2}$}

  \initcurs{240}{0}\SBox\initcurs{232}{8}\SBox\initcurs{232}{16}\SBox
  \initcurs{248}{-8}\SBox\initcurs{256}{-16}\SBox\initcurs{248}{-24}\SBox
  \initcurs{248}{-32}\SBox\initcurs{240}{-40}\SBox
  \put(228,22){$\sliceS{B^s_3}$}
\end{picture}
\end{sideways}
  \caption{\label{D.dttn.zeichnung}
    A minimal deformation shape of $\wild{D}{6}$ having codimension $4$}
\end{figure}

We start with the details of the proof. First we construct $\delta_m$ by
\begin{itemize}
\item setting $\delta_m(U):=1$
\item setting $\delta_m(\trdp{i}{v}):=
  \delta_m(\arsm{\trdp{i}{v}})
  -\delta_m(\trdp{(i-1)}{v})-1$ for any $\trdp{i}{v}$ with
  $\tau^{-1} B^q_j\not\preceq\trdp{i}{v}\not\preceq B^s_j$ for an arbitrary $j$
  or $\trdp{i}{v}\not\preceq B^s_{m-1}$ or $\tau^{-1} B^q_0\not\preceq\trdp{i}{v}$
\item setting $\delta_m(\trdp{i}{w}):=
  \delta_m(\arsm{\trdp{i}{w}})
  -\delta_m(\trdp{(i-1)}{w})-1$ for any $\trdp{i}{w}$ with
  $\tau^{-1} B^q_j\preceq\trdp{i}{w}\preceq B^s_{j+1}$ for an arbitrary $j$
  or $\trdp{i}{w}\not\preceq B^s_{m-1}$ or $\tau^{-1} B^q_0\not\preceq\trdp{i}{w}$
\item making it strictly additive at any other indecomposable
  $\mathcal{Q}$-module $W$ with $\tau^{-1} W\preceq V$
\item setting it equal to zero anywhere else.
\end{itemize}
We illustrate this choice with the case $\wild{D}{6}$ and $m=4$ in
figure~\ref{D.dttn.zeichnung}. We observe that the $\mathcal{Q}$-module $M$
associated to $\delta_m$ has for each $B^q_j$ one direct summand isomorphic to the
$\trdp{i}{v}$ which lies in $B^q_j$ and that for any other direct summand $W$
of $M$ we get $\delta_m(W)=0$. Hence, by lemma \ref{B.deltaFkt}~(c),
$\delta_m$ has codimension $m$.

We define the \emph{defect $d_R(\delta)$ of $\delta$ at $R$} for any
quiver that consists of a tame quiver plus an additional vertex $w$:
For any slice $R$ of $\mathcal{Z}$ and any function $\delta:\mathcal{Z}\to\IZ$
we extend, by abuse of notation, the notion of defect $d_R(\delta)$
of $\delta$ at $R$ from the tame quiver determined by
$\pts{\mathcal{Q}}\setminus\{w\}$ to $\mathcal{Q}$: first we assign to
each vertex $\trdp{i}{j}$, $j\in\pts{\mathcal{Q}}\setminus\{w\}$ in $R$
the weight
\begin{displaymath}
  w_j := \nullr{\mathcal{Q}}{j}
  -\sum_{k\in\pts{\mathcal{Q}}\setminus\{w\}} m(j, k)\nullr{\mathcal{Q}}{k}
\end{displaymath}
where $m(j, k)$ denotes the number of arrows in the Auslander-Reiten quiver
from $\trdp{i}{j}$ to $\trdp{i'}{k}$ with $\trdp{i'}{k}\in R$ for
suitable $i'$ for any pair $(j, k)$. Now we set
\begin{displaymath}
  d_R(\delta) := \sum_{j\in\pts{\mathcal{Q}}\setminus\{w\}}
  w_j \delta_R(j)
\end{displaymath}

We split the proof of the minimality of $\delta_m$ into two
lemmas. The first one applies to all extended tame quivers:

\begin{lem}\label{R.defekt} Let $\mathcal{Q}$ be a quiver of which the
  underlying unoriented graph is of one of the forms
  $\wild{A}{p,q}$, $\wild{D}{n}$, $\wild{E}{6}$,
  $\wild{E}{7}$ or $\wild{E}{8}$, let
  $\delta$ be a deformation shape for $\mathcal{Q}$, let $R$ be a section,
  $p\in\pts{\mathcal{Q}}\setminus\{w\}$ a vertex such that 
  $\tau^{-i}P_p$ for an appropriate $i$ is a source inside $R$.
  \begin{enumerate}
  \item\label{R.defekt.formel} Let $C$ be the generalised Cartan-matrix of
    $\mathcal{Q}$. We have\\
    $\defect{\delta}{R}-\defect{\delta}{\reflect{p}{R}}=$
    $\nullr{\mathcal{Q}}{p}(-\delta_R(p)-\delta_{\reflect{p}{R}}(p)
    -\sum_{q\in\pts{\mathcal{Q}}\setminus\{p,w\}}C_{pq}\delta_R(q))$.
  \item\label{R.defekt.monoton} If $p$ and $w$ are connected by an arrow
    it follows that $\defect{\delta}{R}-\defect{\delta}{\reflect{p}{R}}
    \geq -\nullr{\mathcal{Q}}{p}\delta_R(w)$. Otherwise we have
    $\defect{\delta}{R}-\defect{\delta}{\reflect{p}{R}}\geq 0$
  \item\label{R.defekt.eindeutig} If
    $\defect{\delta}{R}=\defect{\delta}{\reflect{p}{R}}$, then every
    $\delta|_R$ determines $\delta|_{\reflect{p}{R}}$ uniquely and
    vice versa.
  \end{enumerate}
\end{lem}

\noindent {\em Proof.} \hspace*{6pt}
  \begin{enumerate}
  \item
    \begin{eqnarray*}
      & & \defect{\delta}{R} - \defect{\delta}{\reflect{p}{R}} \\
      & = & \sum_{q\in\pts{\mathcal{Q}}\setminus\{p,w\}}
	m(q,p)\nullr{\mathcal{Q}}{p}\delta_R(q) \\
      & & + \left(\nullr{\mathcal{Q}}{p} -
	\sum_{q\in\pts{\mathcal{Q}}\setminus\{w\}} m(p,q)\nullr{\mathcal{Q}}{q}\right)
	  \delta_R(p)
	- \nullr{\mathcal{Q}}{p}\delta_{\reflect{p}{R}}(p) \\
      & = & C_{pp}\nullr{\mathcal{Q}}{p}\delta_R(p)
      + \sum_{q\in\pts{\mathcal{Q}}\setminus\{p,w\}}C_{pq}
          \nullr{\mathcal{Q}}{q}\delta_R(p) \\
      & & + \nullr{\mathcal{Q}}{p}(-\delta_R(p)-\delta_{\reflect{p}{R}}(p)
	-\sum_{q\in\pts{\mathcal{Q}}\setminus\{p,w\}} C_{pq}\delta_R(q))
    \end{eqnarray*}
  \item The subadditivity inequation on $\trdp{i}{p}\in\reflect{p}{R}$
    yields $-\delta_R(p)-\delta_{\reflect{p}{R}}(p)
    -\sum_{q\in\pts{\mathcal{Q}}\setminus\{p,w\}}C_{pq}\delta_R(q)
    -C_{pw}\delta_R(w)\geq 0$.
  \item For a fixed $\delta|_R$, anything in the equation of part (a)
    except $\delta_{\reflect{p}{R}}(p)$ is fixed.
\hspace*{12pt}$\Box$
  \end{enumerate}

We will use the symmetry inherent to the conditions for deformation
shapes to define some notation: Whenever we have chosen an orientation
of a quiver $\mathcal{Q}$ without cycles such that there
is only one sink called $p$ and have a deformation shape
$\delta$ of $\proj{p}$ and $\trdp{(i+1)}{p}$, we define
$k(q):=\max\{k|\trdp{k}{q}\preceq\trdp{i}{p}\}$ for any
$q\in\pts{\mathcal{Q}}$. Then let the \emph{reflection of} $\delta$ be the
deformation shape $\overline{\delta}$ determined by
$\overline{\delta}(\trdp{j}{q}):=\delta(\trdp{(k(q)-j)}{q})$. Given a
slice $R$, the slice $\overline{R}^{\delta}$ defined by
$\overline{R}^{\delta}:=\{\trdp{(k(q)-j)}{q}|\trdp{j}{q}\in R\}$ will be called
the \emph{reflection of} $R$. Clearly $\overline{\delta}$ is a deformation
shape if and only if $\delta$ is. Furthermore, we have
$\defect{\overline{\delta}}{\overline{R}^{\delta}}=-\defect{\delta}{R}$.

Now we can prove the minimality of $\delta_m$ for
$\mathcal{Q}=\wild{D}{n}$:

\begin{lem}\label{dttn} Set $\mathcal{Q}$ equal to $\wild{D}{n}$
  with an orientation fixed by making $z_1$ the only sink and let
  $\delta$ be an arbitrary deformation shape of $U:=\proj{z_1}$
  and $V:=\trdp{((m-1)(c+1)+n-2)}{z_1}$ satisfying
  $\delta\leq\delta_m$. Then the following holds:
  \begin{enumerate}
  \item We have $\defect{\delta}{\sliceQ{B^q_0}}\leq 0$. Moreover,
    $\defect{\delta}{\sliceQ{B^q_0}}=0$ implies
    \mbox{$\delta(\trde{v})=1$}.
  \item If $\defect{\delta}{\reflect{v}{\sliceS{B^s_k}}}\leq -1$
    then $\defect{\delta}{\sliceQ{B^q_k}}\leq 0$. In
    case $\defect{\delta}{\sliceQ{B^q_k}}=0$ we get
    $\delta(\trdp{(k(c+1)+1)}{v})=1$.
  \item We have $\delta=\delta_m$
  \end{enumerate}
\end{lem}

\noindent {\em Proof.} \hspace*{6pt}
  \begin{enumerate}
  \item To obtain an upper bound for the defect we repeatedly add to it the
    subadditivity inequation for a $\preceq$-maximal indecomposable
    $W$ such that $\delta(W)$ occurs in it with multiplicity greater than
    $0$, thus making its multiplicity $0$. We do not replace
    $\delta(\trde{v})$. We use
    the fact that $\delta(W)=0$ for any $W\not\in \supp\delta_m$:
    \begin{eqnarray*}
      \defect{\delta}{\sliceQ{B^q_0}}
      & = & -2\delta(\trde{z_{n-3}}) + \delta(\trde{v}) + \delta(\trde{d})
	+ \delta(\trdp{n-3}{a}) \\
      & & + \delta(\trdp{n-3}{b}) \\
      & \leq & -2\delta(\trde{z_{n-3}}) + \delta(\trde{v}) + \delta(\proj{d})
	+ \delta(\trdp{n-3}{a}) \\
      & & - \delta(\trdp{n-4}{b}) + \delta(\trdp{n-3}{z_1}) \\
      & \leq & \ldots \\
      & \leq & \delta(\trde{v}) - \delta(\proj{z_1}) \leq 0
    \end{eqnarray*}

  \item As in (a) we prove the inequality
    \begin{displaymath}
      \defect{\delta}{\sliceQ{B^q_k}}
      - \defect{\delta}{\reflect{v}{\sliceS{B^s_k}}} \leq
      \delta(\trdp{(k(c+1)+1)}{v})
    \end{displaymath}

  \item We have $\defect{\delta}{\sliceQ{B^q_k}}\leq 0$
    for all $k\in\{0,\ldots,m-2\}$: use that the assumption
    $\defect{\delta}{\reflect{v}{\sliceS{B^s_{(k+1)}}}}=0$
    leads to $\delta(\trdp{(k+1)(c+1)}{v})=\delta(\trdp{(k(c+1)+1)}{v})=1$
    because of the strict additivity of $\delta$ provided by
    lemma~\lemref{R.defekt}{R.defekt.monoton}{b} as well as
    part (a) and (b) of this lemma. This contradicts
    $\delta\leq\delta_m$.

    Furthermore, $\delta_m=\overline{\delta_m}$,
    hence $\overline{\delta}\leq\delta_m$. From
    $\defect{\overline{\delta}}{\sliceQ{B^q_k}}\leq 0$ for
    all $k\in\{0,\ldots,m-2\}$ we obtain
    $\defect{\delta}{\sliceS{B^s_{(k+1)}}}\geq 0$ for all
    $k\in\{0,\ldots,m-2\}$. Thus the fact $\defect{\delta}{\sliceQ{B^q_k}}\leq 0$
    can be sharpened to $\defect{\delta}{\sliceQ{B^q_k}}=0$. This shows
    $\delta(\trdp{(k(c+1)+1)}{v})=1$ for all $k\in\{0,\ldots,m-2\}$.
    But from $\defect{\overline{\delta}}{
    \reflect{v}{\sliceS{B^s_{(k+1)}}}}\leq -1$ we obtain
    $\defect{\delta}{\reflectInv{v}{\sliceQ{B^q_k}}}\geq 1$ for
    all $k\in\{0,\ldots,m-2\}$. Thus the subadditivity inequality is
    in each $\trdp{(k(c+1)+1)}{v}$ proper subadditive while
    $\delta(\trdp{(k(c+1)+1)}{v})=1$. This forces
    $\delta$ to have codimension at least $m$. But any $\delta < \delta_m$
    must have codimension strictly smaller than $m$. We conclude
    $\delta=\delta_m$.
\hspace*{12pt}$\Box$

  \end{enumerate}

We sketch for the other quivers that consist of a tame quiver as a full
subquiver plus an additional vertex $w$ the construction of the families
of minimal disjoint degenerations of arbitrary codimension:
In the case of $\wild{E}{6}$ set $U:=\proj{z}$, $V:=\trdp{(8m-4)}{z}$,
$B^q_i := \trdp{(8i+2)}{z}$, $B^s_i := \trdp{(8i+1)}{z}$ and
construct $\delta_m$ in the same way as in the case
$\wild{D}{n}$. For $\wild{E}{7}$ set $U:=\proj{z}$, $V:=\trdp{(15m-9)}{z}$,
$B^q_i := \trdp{(15i+3)}{z}$, $B^s_i := \trdp{(15i+2)}{z}$. The quiver
$\wild{E}{8}$ requires some, the quivers $\wild{A}{p,q}$ a lot of
further technical adaption and can be found in \cite{Olbricht:2004}.

The deformation shapes of a quiver with underlying unoriented graph $S$
have a nicer property than those of $\wild{D}{n}$: We fix the orientation
of $S$ such that $z$ is the only sink. To a pair $U:=\proj{z}$,
$V:=\trdp{(2m+1)}{z}$ we can construct minimal deformation shapes of
$U$ and $V$ to any codimension $n\in\{1,\ldots,m\}$. Let $\delta_{m,n}$
be additive in any indecomposable except
\begin{itemize}
\item $\delta_{m,n}(U)=\delta_{m,n}(\tau V)=1$.
\item $\delta_{m,n}(\trdp{(4i+1)}{d})
  =\delta_{m,n}(\trdp{(4i+2)}{d})=1$ for
  $i\in\{i\in\IN|2i+2\leq n\}$ and, if $n$ is even,
  $\delta_{m,n}(\trdp{i}{d})=1$ for $i\in\{2n-1,\ldots,2m-1\}$.
\item $\delta_{m,n}(\trdp{(4i+3)}{e})
  =\delta_{m,n}(\trdp{(4i+4)}{e})=1$ for
  $i\in\{i\in\IN|2i+3\leq n\}$ and, if $n$ is odd,
  $\delta_{m,n}(\trdp{i}{e})=1$ for $i\in\{2n-1,\ldots,2m-1\}$.
\item $\delta_{m,n}$ is zero for any other indecomposable in the
  $\tau$-orbit of $\proj{d}$ or $\proj{e}$.
\end{itemize}

Then $\delta_{m,n}$ is a deformation shape of $U$ and $V$
and has codimension $m$. We can prove the minimality by a lemma analogous
to lemma~\ref{dttn}.

\newpage

\bibliography{zitate}

\end{document}